\documentclass[reqno,10pt]{amsart}
\usepackage{url}
\usepackage{amssymb,amsthm,amsmath,amsfonts,amstext,mathrsfs}
\usepackage{latexsym}
\usepackage[all]{xy}
\newtheorem{thm}{Theorem}
\newtheorem*{lem}{Lemma}

\input xy
\xyoption{all}

\def\d{\,{\rm{d}}}
\def\m{\,{\mathfrak{m}}}

\title[Giedrius Alkauskas]
{The Minkowski $?(x)$ function and Salem's problem\\
La fonction $?(x)$ de Minkowski et probl\`{e}me de Salem}
\author[Giedrius Alkauskas]{Giedrius Alkauskas}
\begin{document}

\begin{abstract} R. Salem (Trans. Amer. Math. Soc. 53 (3) (1943) 427-439) asked whether the Fourier-Stieltjes transform of the Minkowski question mark function $?(x)$ vanishes at infinity. In this note we present several possible approaches towards the solution. For example, we show that this transform satisfies integral and discrete functional equations. Thus, we expect the affirmative answer to Salem's problem. In the end of this note we show that recent attempt to settle this question (S. Yakubovich, C. R. Acad. Sci. Paris, Ser. I {\bf 349} (11-12) (2011) 633-636) is fallacious.\\

\noindent {\sc R\'{e}sum\'{e}.} R. Salem (Trans. Amer. Math. Soc. 53 (3) (1943) 427-439) demande si la transform\'{e}e de Fourier-Stieltjes  de la fonction point d'interrogation de  Minkowski $?(x)$ s'annule \`{a} l'infini. Dans cette note nous pr\'{e}sentons plusieurs approches afin de r\'{e}soudre cette question. Nous montrons par exemple que cette transforme\'{e}e satisfait des \'{e}quations fonctionnelles discr\`{e}tes et enti\`{e}res. Ainsi, nous conjecturons une r\'{e}ponse positive au probl\`{e}me de Salem. A la fin de cette note, nous montrons qu'une  tentative r\'{e}cente pour r\'{e}pondre \`{a} cette question (S. Yakubovich, C. R. Acad. Sci. Paris, Ser. I {\bf 349} (11-12) (2011) 633-636) est en fait incorrecte.
\end{abstract}

\maketitle
\begin{center}
Mathematical Analysis/Number Theory
\end{center}

\section{Salem's problem}
The Minkowski question mark function $?(x):[0,1]\mapsto[0,1]$ is defined by
\begin{eqnarray*}
?([0,a_{1},a_{2},a_{3},\ldots])=2\sum\limits_{i=1}^{\infty}(-1)^{i+1}2^{-\sum_{j=1}^{i}a_{j}},\quad a_{i}\in\mathbb{N};
\end{eqnarray*}
$x=[0,a_{1},a_{2},a_{3},\ldots]$ stands for the representation of $x$ by a (regular) continued fraction. The function $?(x)$ is continuous,
strictly increasing, and singular. The extended Minkowski question mark function is defined by
$F(x)=?(\frac{x}{x+1})$, $x\in[0,\infty)$. Thus, for $x\in[0,1]$, we
have $?(x)=2F(x)$. The function $F(x)$ satisfies functional
equations
\begin{eqnarray}
2F(x)=\left\{\begin{array}{c@{\qquad}l} F(x-1)+1 & \mbox{if}\quad x\geq 1,
\\ F({x\over 1-x}) & \mbox{if}\quad 0\leq x<1. \end{array}\right.\label{distr}
\end{eqnarray}
This implies $F(x)+F(1/x)=1$. As was proved by Salem \cite{Salem}, the function $?(x)$ satisfies H\"{o}lder condition of order
$\alpha=(2\log \frac{\sqrt{5}+1}{2})^{-1}\log 2=0.7202_{+}$.
The Laplace-Stieltjes transform of $?(x)$ is defined by \cite{ga}
\begin{eqnarray*}
\m(t)=\int\limits_{0}^{1}e^{xt}\d ?(x),\quad t\in\mathbb{C}.
\end{eqnarray*}
This is an entire function. The symmetry property $?(x)+?(1-x)=1$ implies $\m(t)=e^{t}\m(-t)$.
Let $d_{n}=\m(2\pi i n)$, $n\in\mathbb{N}$. Because of the symmetry property we have $d_{n}\in\mathbb{R}$, and thus
\begin{eqnarray}
d_{n}=\int\limits_{0}^{1}\cos(2\pi n x)\d ?(x),\quad n\in\mathbb{N}.\label{four}
\end{eqnarray}
In 1943 Rapha\"{e}l Salem \cite{Salem} posed the following problem: prove or disprove that  $d_{n}\rightarrow 0$, as $n\rightarrow\infty$.
The question to determine whether Fourier transform of a given measure vanishes at infinity is a very delicate question whose answer depends on an intrinsic structure of this measure. There are various examples for both cases of behaviour \cite{zygmund}.  As was noted in \cite{Salem}, the general
theorem of Wiener \cite{zygmund} about Fourier-Stieltjes
coefficients of continuous monotone functions with known modulus of continuity and the Cauchy-Schwartz inequality imply that $\sum_{n=1}^{N}|d_{n}|=O(N^{1-\alpha/2})$. Thus, $|d_{n}|\ll
n^{-0.3601}$ on average. Via a partial summation and standard calculations we get that
\begin{eqnarray}
\sum\limits_{n=1}^{\infty}\frac{d_{n}}{n^{\sigma}}\text{ converges absolutely for }\sigma>1-\frac{\alpha}{2}=0.6398_{+},\text{ and }
?(x)-x=\sum\limits_{n=1}^{\infty}\frac{d_{n}}{\pi n}\cdot\sin(2\pi nx),\quad x\in[0,1]. \label{conv}
\end{eqnarray}
Note that $|\sum_{n=1}^{N}d_{n}|\leq 2\int_{0}^{1}|e^{2\pi i x}-1|^{-1}\d ?(x)$ which is finite, since $?(x)=1-?(1-x)\asymp 2^{-1/x}$ as $x\rightarrow 0_{+}$. Thus, we inherit that the Dirichlet series $\sum_{n=1}^{\infty}d_{n}n^{-\sigma}$ converges (conditionally) for $\sigma>0$. \\
\indent The purpose of this note is to disseminate the knowledge of Salem's problem to a wider audience of mathematicians. We contribute to
this topic with two new results. Vaguely speaking, they show that the coefficients $d_{n}$ behave in the same manner as they behave ``on average";
hence the answer to Salem's problem most likely is positive. Let, as usual, $J_{\nu}(\star)$ stand for the Bessel function with index $\nu$.
\begin{thm} {\rm (Integral functional equation)}. The function $\m(it)$ satisfies the following identity:
\begin{eqnarray*}
\frac{i\m(is)}{2e^{2is}-e^{is}}&=&\int\limits_{0}^{\infty}\m'(it)J_{0}(2\sqrt{st})\d
t,\quad s>0.
\end{eqnarray*}
The integral is conditionally convergent. \label{thm1}
\end{thm}
Note that $\m(t)$ satisfies analogous integral equation
on the negative real line \cite{ga}. Theorem \ref{thm1}, however, cannot be deduced from the latter by standard methods. If we formally
pass to the limit $s\rightarrow\infty$ under the integral, the bound $|J_{0}(2\sqrt{st})|\ll (st)^{-1/4}$ would imply $\m(is)\rightarrow 0$.
 The same conclusion follows if we formally take the limit $s\rightarrow 0$. Unfortunately, this conditionally convergent integral cannot
 be dealt this way. In fact, let $\mathfrak{n}(it)=e^{it}$. Then (\ref{bess}) shows that
$i\mathfrak{n}(is)e^{-2is}=\int_{0}^{\infty}\mathfrak{n}'(it)J_{0}(2\sqrt{st})\d t$ for $s>0$. Now the formal passage to the limit 
$s\rightarrow\infty$ gives the false result $e^{is}\rightarrow 0$. Therefore, if
the solution of Salem's problem based on Theorem \ref{thm1} is found, it should deal with the factor $(2e^{2is}-e^{is})^{-1}$, as opposed to $e^{-2is}$.
The behaviour of $?(x)$ at $x=0$ and $x=1$ is of importance as well. Theorem \ref{thm1} has a discrete analogue.
\begin{thm}{\rm (Discrete functional equation)}. For any $m\in\mathbb{N}$ we have the following identity:
\begin{eqnarray*}
d_{m}=\int\limits_{0}^{1}\cos\Big{(}\frac{2\pi m}{x}\Big{)}\d x+
2\sum\limits_{n=1}^{\infty}d_{n}\cdot\int\limits_{0}^{1}\cos(2\pi n x)\cos\Big{(}\frac{2\pi m}{x}\Big{)}\d x.
\end{eqnarray*}
This sum is majorized by the series $Cm\sum_{n=1}^{\infty}|d_{n}|n^{-3/4}$ (see (\ref{conv})) with an absolute constant $C$.\label{thm2}
\end{thm}
The theorem of Salem and Zygmund \cite{zygmund} shows that $d_{n}=o(1)$ implies $\m(it)=o(1)$; this is a general fact for the
Fourier-Stieltjes transforms of non-decreasing functions. Another idea how to tackle Salem's problem is to approach
it via the above system of infinite linear identities. This demands a detailed study of the integral $P(a,b)=\int_{0}^{1}\cos(a/x+bx)\d x$.
Its exact asymptotics can be given in terms of elementary functions if $b>(1+\epsilon)a$, or $b<(1-\epsilon)a$ (b can be negative), or $b=a$,
where $\epsilon>0$ is fixed.
The transition area $b\sim a$ exhibits a more complex behaviour. One can nevertheless give exact asymptotics in terms of Fresnel sine and cosine integrals, and this asymptotics is also valid in the transition area. These investigations are due to N. Temme \cite{Temme}. Possibly, the full strength of these results
can solve Salem's problem; our joint project with N. Temme is in progress.
\section{The proofs}
\noindent{\bf Proof of Theorem \ref{thm1}}. First, we will show that the integral does converge relatively. Indeed, let $X>0$. Then
\begin{eqnarray*}
A(s,X):=\int\limits_{0}^{X}\mathfrak{m}'(it)J_{0}(2\sqrt{st})\d
t=-i\int\limits_{0}^{X}J_{0}(2\sqrt{st})\d\mathfrak{m}(it)=-iJ_{0}(2\sqrt{sX})\m(iX)+i-
i\int\limits_{0}^{X}\mathfrak{m}(it)J_{1}(2\sqrt{st})\frac{s^{1/2}}{t^{1/2}}\d t.
\end{eqnarray*}
Let
\begin{eqnarray}
\widehat{\mathfrak{m}}(T)=\int\limits_{0}^{T}\mathfrak{m}(it)\d
t=\int\limits_{0}^{1}\frac{e^{ixT}-1}{ix}\d ?(x).\text{ This implies }
|\widehat{\mathfrak{m}}(T)|\leq2\int\limits_{0}^{1}x^{-1}\d
?(x)=5,\text{ for }T\geq 0.
\label{antii}
\end{eqnarray}
We can continue:
\begin{eqnarray}
A(s,X)=-iJ_{0}(2\sqrt{sX})\m(iX)+i-i\int\limits_{0}^{X}J_{1}(2\sqrt{st})\cdot\frac{s^{1/2}}{t^{1/2}}\d\widehat{\mathfrak{m}}(t)\nonumber\\=
-iJ_{0}(2\sqrt{sX})\m(iX)+i-iJ_{1}(2\sqrt{sX})\frac{s^{1/2}}{X^{1/2}}\widehat{\m}(X)
+\frac{i}{2}\int\limits_{0}^{X}\widehat{\mathfrak{m}}(t)\Big{(}J_{0}(2\sqrt{st})\frac{s}{t}-
J_{2}(2\sqrt{st})\frac{s}{t}-J_{1}(2\sqrt{st})\frac{s^{1/2}}{t^{3/2}}\Big{)}\d
t.\label{relative}
\end{eqnarray}
The function under integral is bounded in the neighborhood of $t=0$ since
$\widehat{\mathfrak{m}}(t)$ has a first order zero at $t=0$, and for
$\nu\in\mathbb{N}_{0}$, $J_{\nu}(u)$ has a zero of order $\nu$ at
$u=0$ (thus, no zero for $\nu=0$). Further, we have the bound for the Bessel function
$|J_{\nu}(u)|\ll u^{-1/2}$ as $u\rightarrow\infty$, $\nu$ is fixed.
Thus, the function under integral is $\ll t^{-5/4}$
for $t>1$, hence the integral (\ref{relative}) converges absolutely. Therefore there
exists a finite limit $A(s)=\lim\limits_{X\rightarrow\infty}A(s,X)$, and the integral in Theorem \ref{thm1} converges
conditionally. In fact, we used only the properties (\ref{antii})
and $|\m(it)|\leq 1$. Further, we take $X=\infty$ and substitute
(\ref{antii}) into (\ref{relative}). We get
\begin{eqnarray}
A(s)=
i+\frac{i}{2}\int\limits_{0}^{\infty}\Big{[}\int\limits_{0}^{1}\frac{e^{ixt}-1}{ix}\d
?(x)\Big{]}\Big{(}J_{0}(2\sqrt{st})\frac{s}{t}-
J_{2}(2\sqrt{st})\frac{s}{t}-J_{1}(2\sqrt{st})\frac{s^{1/2}}{t^{3/2}}\Big{)}\d
t.
\label{inter}
\end{eqnarray}
The double integral converges absolutely. Indeed,
$(e^{ixt}-1)(ix)^{-1}=\int_{0}^{t}e^{ixu}\d u\Rightarrow \big{|}(e^{ixt}-1)(ix)^{-1}\big{|}\leq \min\big{\{}t,2/x\big{\}}$.
Now we easily obtain an absolute convergence: just use the bound $t$
for $t\leq 1$ and the bound $2/x$ for $t\geq 1$. So Fubini's theorem
allows us to interchange the order of integration in (\ref{inter}).
After going backwards by integrating by parts, we get
\begin{eqnarray*}
A(s)=
i+\frac{i}{2}\int\limits_{0}^{1}\Big{[}\int\limits_{0}^{\infty}(e^{ixt}-1)
\Big{(}J_{0}(2\sqrt{st})\frac{s}{t}-
J_{2}(2\sqrt{st})\frac{s}{t}-J_{1}(2\sqrt{st})\frac{s^{1/2}}{t^{3/2}}\Big{)}\d
t\Big{]}\frac{\d ?(x)}{ix}=\int\limits_{0}^{1}\Big{[}\int\limits_{0}^{\infty}x
e^{ixt}J_{0}(2\sqrt{st})\d t\Big{]}\d ?(x).
\end{eqnarray*}
For $x,s>0$, we have the classical integral \cite{watson}
\begin{eqnarray}
x\int\limits_{0}^{\infty}e^{ixt} J_{0}(2\sqrt{st})\d
t=i e^{-\frac{is}{x}}.\label{bess}
\end{eqnarray}
So, using $F(x)+F(1/x)=1$, functional equations (\ref{distr}), and the symmetry property for $\m(t)$, we obtain
\begin{eqnarray*}
A(s)=2i\int\limits_{0}^{1}e^{-\frac{is}{x}}\d
F(x)\mathop{=}^{x\mapsto\frac{1}{x}}2i\int\limits_{1}^{\infty}
e^{-isx}\d F(x)
\mathop{=}^{(\ref{distr})}2i\int\limits_{0}^{1}\sum\limits_{n=1}^{\infty}\frac{e^{-is(x+n)}}{2^{n}}\d
F(x)=\frac{i\mathfrak{m}(-is)}{2e^{is}-1}=\frac{i\mathfrak{m}(is)}{2e^{2is}-e^{is}}.\quad\square
\end{eqnarray*}

\noindent{\bf Proof of Theorem \ref{thm2}}. We know that $?(x)-x$ can be expressed by the absolutely uniformly
convergent series (\ref{conv}), which can be integrated term-by-term. Let, for $0<\epsilon<1$,
\begin{eqnarray*}
\widehat{A}(s,\epsilon)=i\int\limits_{\epsilon}^{1}e^{-\frac{is}{x}}\d
?(x)=2i\int\limits_{1}^{1/\epsilon}e^{-isx}\d F(x).
\end{eqnarray*}
We know that
$\lim\limits_{\epsilon\rightarrow 0_{+}}\widehat{A}(s,\epsilon)=A(s)$.
Thus,
\begin{eqnarray}
\widehat{A}(s,\epsilon)&=&i\int\limits_{\epsilon}^{1}e^{-\frac{is}{x}}\d[(?(x)-x)+x]=i\int\limits_{\epsilon}^{1}e^{-\frac{is}{x}}\d x
-ie^{-\frac{is}{\epsilon}}(?(\epsilon)-\epsilon)+\int\limits_{\epsilon}^{1}(?(x)-x)e^{-\frac{is}{x}}\cdot\frac{s\d x}{x^{2}}\nonumber\\
&=& i\int\limits_{\epsilon}^{1}e^{-\frac{is}{x}}\d x
-ie^{-\frac{is}{\epsilon}}(?(\epsilon)-\epsilon)+\sum\limits_{n=1}^{\infty}\frac{d_{n}}{\pi n}\int\limits_{\epsilon}^{1}\sin(2\pi nx)e^{-\frac{is}{x}}
\cdot\frac{s\d x}{x^{2}}\nonumber\\
&=&i\int\limits_{\epsilon}^{1}e^{-\frac{is}{x}}\d
x-ie^{-\frac{is}{\epsilon}}(?(\epsilon)-\epsilon)
+i\sum\limits_{n=1}^{\infty}\frac{d_{n}}{\pi n}\sin(2\pi n
\epsilon)e^{-\frac{is}{\epsilon}}+2i\sum\limits_{n=1}^{\infty}d_{n}\int\limits_{\epsilon}^{1}\cos(2\pi
nx)e^{-\frac{is}{x}}\d x.\label{ser}
\end{eqnarray}
Take the imaginary part. Two series converge absolutely and uniformly with respect to $\epsilon$. This follows from (\ref{conv}) and
\begin{lem}
Let $a>0$. For a certain absolute constant $C>1$ and any $\epsilon\in[0,1]$ we have the following bound:
\begin{eqnarray*}
\Big{|}\int\limits_{\epsilon}^{1}\cos\Big{(}bx+\frac{a}{x}\Big{)}\d
x\Big{|}<\left\{\begin{array}{c@{\qquad}l}C(a+1)b^{-3/4} & \mbox{if}\quad b\geq  2\pi ,
\\C(a+1)|b|^{-1} & \mbox{if}\quad b\leq -2\pi. \end{array}\right.
\end{eqnarray*}
\end{lem}
\noindent If $2\pi\leq b\leq (C(a+1))^{4/3}$, the first bound is trivial since the integral is $\leq 1$. If $b>(C(a+1))^{4/3}$, we have the case 
$b>(1+\epsilon)a$ in the
aforementioned cosine integral $P(a,b)$. In this case the asymptotics has only contributions from a stationary saddle point $x_{0}=(a/b)^{1/2}$
(if the latter is in the neighborhood of $\epsilon$) and the end point $x_{1}=1$. We deal with the case $b<0$ similarly. Thus, the above bounds follow
from results and techniques in \cite{Temme}; the exponent $3/4$ is the best possible and cannot be increased. This proves the Lemma. Note that $2\cos(a/x)\cos(bx)=\cos(a/x+bx)+\cos(a/x-bx)$.
Therefore for fixed $s$, the last series in (\ref{ser}) is majorized by $\hat{C}(s+1)\sum_{n=1}^{\infty}|d_{n}|n^{-3/4}$, and one can pass to the limit $\epsilon\rightarrow 0_{+}$ in (\ref{ser}) elementwise. This yields
\begin{eqnarray*}
\Re\frac{\m(is)}{2e^{2is}-e^{is}}=\Im A(s)=\lim\limits_{\epsilon\rightarrow 0_{+}}\Im\widehat{A}(s,\epsilon)=\int\limits_{0}^{1}\cos\Big{(}\frac{s}{x}\Big{)}\d
x+2\sum\limits_{n=1}^{\infty}d_{n}\cdot\int\limits_{0}^{1}\cos(2\pi
nx)\cos\Big{(}\frac{s}{x}\Big{)}\d x.
\end{eqnarray*}
Now we finish with the substitution $s=2\pi m$, $m\in\mathbb{N}$. $\square$
\appendix\section{}
The proof of Salem's problem presented in \cite{yak} is fallacious and cannot be fixed. Indeed, if we track down what properties of $?(x)$ are
used in the proof, we find that the author only uses the fact that $?(x)$ is of bounded variation,
that $?(0)=0$, $?(1)=1$, $\int_{0}^{1}?(x)x^{-5/4}\d x<+\infty$, and $?(x)=2F(x)=2+O(x^{-3})$ as $x\rightarrow+\infty$.
Moreover, the last property is not needed as well, as we will now explain. The author writes
$\int_{0}^{1}=\int_{0}^{\infty}-\int_{1}^{\infty}$ (formula (18) in \cite{yak}), and derives asymptotic expansions of both integrals on the
right side. But our main concern is still the asymptotics for the integral $\int_{0}^{1}f(x)D(\tau,x)\d x$, formula \cite{yak}, (17). Now it is
obvious that we can extend the function $f(x)$, initially defined for $[0,1]$, to the interval $[1,\infty)$
in an almost arbitrary way. Eventually, if the asymptotic results we obtain are correct, two contributions from $[1,\infty)$
will annihilate one another. Nevertheless, the properties which were really used are obviously insufficient for the Fourier-Stieltjes
transform of a singular measure to vanish at infinity; the Cantor ``middle-third" distribution is a counterexample.\\
\indent Now we will indicate where the mistake is. Assume that D. Naylor's result, which is the basis of authors results, is true
(\cite{yak}, formulas (8)-(9)). This would give the following consequence. Let $f(x)$ be continuous, $f(0)=0$, $f(x)=O(x^{b})$, as
$x\rightarrow 0_{+}$ for a certain $b>0$, and $f(x)=0$ for $x\geq 1$. Then for a fixed $f$,
\begin{eqnarray}
\int\limits_{0}^{1}K_{i\tau}(x)f(x)\frac{\d x}{x}=O(e^{-\pi\tau/2}\tau^{-N})\text{ for every }N\geq 1,\text{ as }\tau\rightarrow\infty.
\label{ineqq}
\end{eqnarray}
In fact, for $\tau>0$ and any $x\in(0,1]$ we have $K_{i\tau}(x)\sim-(2\pi)^{-1/2}e^{-\pi\tau/2}\tau^{-1/2}
\sin\big{[}\tau\log(\frac{ex}{2\tau})\big{]}$.
One can extract further asymptotic terms which are of order $\tau^{-1}$, $\tau^{-2}$ etc. smaller than the first one. Thus, direct calculation
shows that the estimate
(\ref{ineqq}) cannot hold for every such function $f(x)$ and $N=2$ (just consider functions $f$ such that $f(x)=0$ for $x\in[0,1/2]$).
And indeed, Naylor's expansion works for
$0<\beta<\frac{\pi}{2}$. Meanwhile the case $\beta=\frac{\pi}{2}$, which is essential for the argument of \cite{yak}
to work, is not allowed, and, moreover, this expansion in case $\beta=\frac{\pi}{2}$ is false, as we have just seen.
\section*{Acknowledgements}\noindent The author sincerely thanks Nico Temme for the help with the Lemma. This work was supported
by the Lithuanian Science Council whose postdoctoral fellowship 
is being funded by European Union Structural Funds project ``Postdoctoral Fellowship Implementation in Lithuania".

\noindent Vilnius University, Department of Mathematics and Informatics, Naugarduko 24, LT-03225 Vilnius, Lithuania.
{\tt giedrius.alkauskas@gmail.com}

\end{document}